\def\hod{\mbox{HOD}}
\def\ZFC{\mbox{ZFC}}
\def\PD{\mbox{PD}}

\def\ran{\rm ran}

\def\phi{\varphi}\def\mm{\mathcal{M}}
\newcommand{\open}{\Bbb}

\newcommand{\oN}{{\open N}}

\newcommand{\oR}{{\open R}}

\def\mm{\mathfrak{M}}

\def\mn{\mathfrak{N}}

\def\={=\!\!}

\documentclass[12pt]{article}
\usepackage{amsthm,hyperref,comment} 
\usepackage{graphicx,enumitem} 
\usepackage{amssymb,latexsym,times,amsmath}
\newcommand{\psfrag}[2]{}
\usepackage{makeidx,proof}
\usepackage{xcolor}
\theoremstyle{plain}

  \theoremstyle{definition}

  \theoremstyle{remark}

\begin{document}

\title{Second order logic and Set Theory Redux}
\author{Jouko V\"a\"an\"anen\thanks{I have received support from the Academy of Finland, decision number 368671, and from the European Research Council (ERC), grant 101020762.}}
\maketitle


\begin{abstract}
We argue that second order logic is a weaker form of set theory, despite the fact that the former is formalized in a second order langauge and the latter in  a first order language. Along the way, we review the history of interactions between second order logic and set theory, as well as some modern trends.
\end{abstract}

\section{Introduction}
Second order logic is a fragment of set theory. It is powerful in what it can express,  but weaker than set theory. This is clear if we consider informal second order logic and informal  set theory, and compare them. On this informal level it is obvious that second order logic is weaker than set theory. After all, second order logic quantifies over subsets of the domain while set theory quantifies over subsets, sets of subsets, sets of sets of subsets, etc, ad infinitum. To make the comparison more precise, we may consider formalizations.  Formal second order logic with the Comprehension Schemas as axioms and Henkin models as semantics is very much like formal set theory with $\ZFC$ as axioms and binary structures $(M,E)$ as semantics. The biggest difference is that set theory has the  Replacement Schema and second order logic has nothing comparable to it. If we focus on \emph{semantic} consequence and take full Henkin models as the semantics of second order logic, it is comparable to set theory with the structures $(V_\alpha,\in)$ as models. 

 This paper argues that second order logic is nothing more and nothing less than  a powerful fragment of set theory. 
 In set theory we can talk about endless iterations of the power-set operation while in second order logic we  apparently have limit ourselves to just one application of the power-set operation. The statement that second order logic is a fragment of set theory cannot be \emph{literally} true because the formal language of second order logic is different from the formal language of set theory. How can we compare two foundational theories which are given in totally different formal languages? If they were presented in the same formal language they could be readily compared.  This is a basic difficulty in comparing the two---second order logic and set theory---and leads to constant misunderstandings and misrepresentations. This paper aims at clarifying the situation.

A further complication in comparing second order logic and set theory is the following: The basic idea of set theory is that there is one underlying concept, the concept of set, and everything is built from this. Numbers, relations, functions, etc are all sets. The point is \emph{not} that everything must be written in the language of set theory but rather that everything 
\emph{can} be so written, if desired. This ``monist" foundation is sometimes considered ``wrong" and generating false connections between mathematical objects with no connections to mathematical practice\footnote{See e.g. \cite{Simons2005-SIMAST}}. On the other hand, there are undoubtedly also benefits from  everything being built from the same elements, for example, results can be easily compared. If a result e.g. in group theory or topology is based on a principle that cannot be recognized even after considerable effort as true, one may ask whether the principle uses in an unavoidable fashion some principle of set theory, e.g. the Continuum Hypothesis or the existence of a measurable cardinal. Set theory provides a universal language for mathematics. Whatever is proved in this universal language can---ideally---be generally accepted as true, and whatever is unprovable in set theory is probably unprovable in a more restricted mathematical theory as well.

In contrast to the monism of set theory, second order logic has a structuralist or pluralist flavor in the following sense: From the point of view of second order logic it is natural to think that the foundations of mathematics should not be just one big ``theory of everything" but rather a myriad of small pieces of foundations perfectly designed for different purposes. The basic textbook \cite{zbMATH00050939} of second order logic  has the title ``Foundations without Foundationalism" and argues strongly against the idea that the right foundation of mathematics is  axiomatic  set theory formulated in first order logic. It is argued in \cite{zbMATH00050939} that classical mathematics is an arena where second order logic and pluralism reign. However, mathematics has a great degree of unity even if everything is not reduced to the epsilon relation. Mathematics does not consist of isolated islands. It is not clear how second order logic accounts for that.

The conclusion of our preliminary comparison of set theory and second order logic is that set theory takes the extreme position of everything living in a huge set theoretical universe while second order logic takes another kind of extreme position: everything has their own world and there is no single all-encompassing universe. A compromise is to adopt the set theory position and use second order logic as a measure of definability in order to delineate the set theoretic universe. A good example of this is the rich theory of the structural properties of second order (over the natural numbers) definable sets of real numbers in the universe of set theory with infinitely many Woodin cardinals. Here the large cardinals of set theory are used to decide otherwise undecidable classical questions, such as Lebesgue measurability and prefect set properties, about (projective) sets of real numbers. A different kind of compromise is to adopt the second order logic position on a large domain---so large that it covers all the relevant structures one may need. An example might be $V_\kappa$ where $\kappa$ is the smallest inaccessible cardinal. This structure has a second order characterization i.e. there is a second order sentences with this structure, up to isomorphism, as its only model.

Let us return to comparing second order logic and set theory, but first  on the totally informal level. The informal level is good because it 
is---by default---free from the effects that differences of formal language cause. The informal level is, of course,  beset with paradoxes.  Nevertheless, it gives a useful  intuitive first picture of the situation. One cannot formalize without an informal idea of what  it is that is being formalized.  As is well-known, formalization does not necessarily free us from the paradoxes, apart from the most obvious ones, but it is the best we can do. When Zermelo first axiomatized set theory  \cite{zbMATH02639893} his explicit goal was to rise above the paradoxes.

On the informal level there is so-called naive set theory, meaning set theory without formalization, the received starting point of becoming acquainted with the set-theoretical view of mathematics. An example is Cantor's original approach to set theory and later e.g.  \cite{MR114756}. We can similarly conceive of naive second order logic, or informal second order logic. This is the way things are usually presented, because second order logic, as its name suggests, is a logic, hence a formal language. However, we can easily conceive of naive second order \emph{mathematics}, if not second order \emph{logic}.

\emph{Naive set theory} has one single type of objects, namely sets. Sets are (i.e. the domain of sets is)  closed under such operations as indexed unions and intersections as well as cartesian products. Sets are (i.e. the domain of sets is) closed under powersets, meaning that if we have a set $A$ we can collect all subsets  of $A$   (that are in our domain) to one new set called the powerset of $A$.    Furthermore we can use any well-defined property to cut from any given set the elements that satisfy the property and form a new set from such elements. But now the weakness of the informal approach emerges: it is not clear what ``well-defined property" means. The vagueness of the concept of ``well-defined" can be clarified, following Skolem \cite{Skolem1923}, by taking a step towards formalized set theory, by stipulating that well-defined properties at least contain atomic (i.e. $x=y$ and $ x\in y$) properties and are closed under Boolean operations and existential and universal quantifiers (over sets). This will suffice for the needs of developing mathematics in set theory. Normally we do not refer to a ``domain" in naive set theory. In a sense, in naive set theory, everything belongs to the ``domain" and the ``domain" itself is not a set. This implies that the idea of a universal set is not part of naive set theory\footnote{``Nothing contains everything" (\cite{MR114756})}. Not only because it leads immediately to paradoxes, but also because it is not needed. The set of natural numbers can be defined as the intersection of all sets which contain the empty set and are closed under the function which maps $x$ to $x\cup\{x\}$. Naively thinking, such sets exist, because $\emptyset, \{\emptyset\},\{\emptyset,\{\emptyset\}\},\ldots$ is one such.  The set of real numbers can be built from the powerset of the set of natural numbers. One can develop a lot of mathematics in set theory on this naive level (see e.g. \cite{MR114756}). 

\emph{Naive second order mathematics} has  objects we call individuals, objects we call sets (or classes)  of individuals, for each natural number $n$  objects we call $n$-ary relations, and similarly for $n$-ary functions. With any well-defined property of individuals we can introduce the set (or class) of individuals with that property. What ``well-defined" means is left vague as in naive set theory. We can, taking a step towards formalized second order logic, stipulate that well-defined properties contain at least  atomic  properties and are closed under Boolean operations and existential and universal quantifiers over individuals, sets, relations and functions. Similarly we can use well-defined properties of $n$-tuples of  individuals to define $n$-ary relations, and functional properties to define functions.  This will suffice for most purposes. 

The set of natural numbers can be defined in naive second order mathematics as a structure  $N$ with a unary (``successor") function $s$ and a distinguished element (``zero") $0$ satisfying  the familiar second order induction axiom. It is easy to argue informally that such a structure is unique up to isomorphism in the following sense: If  among the sets (i.e. classes) that we have there are two triples $(N,s,0)$ and $(N',s',0')$ which both satisfy the said definition of the natural numbers, then there is an isomorphism $$f:(N,s,0)\cong(N',s',0').$$ The function $f$ can be defined as (or from) the intersection of all binary relations between $N$ and $N'$ satisfying some obvious closure properties. This uniqueness is a phenomenon, characteristic of second order mathematics, called categoricity. The structure of the natural numbers is said to be defined in a categorical way because it is defined up to isomorphism (only). 

The set of real numbers can be similarly  defined as the unique, up to isomorphism, ordered field $(R,+,\times,0,1,<)$ which satisfies the second order property that every non-empty bounded subset has a least upper bound. If  among the sets (i.e. classes) that we have there are two 6-tuples $(R,+,\times,0,1,<)$ and $(R',+',\times',0',1',<')$ which both satisfy the said definition of the real numbers, then, famously (\cite{zbMATH02659709}), there is an isomorphism $$f:(R,+,\times,0,1,<)\cong(R',+',\times',0',1',<').$$
One can develop a lot of mathematics in a second order way on this naive level along these lines. In particular, relevant structures are characterized categorically in naive second order mathematics.

There is no problem in extending the idea of second order mathematics to third order mathematics, and even further. This leads eventually to simple type theory. This corresponds quite closely to how classical mathematics works.

Let us now compare the naive versions of set theory and second order mathematics. In both cases we have some objects and we have ideas how to build new objects. In naive set theory the objects are of one kind only, viz. sets. In naive second order mathematics the objects are of several types, viz. individuals, sets, relations, and functions, of different arities.  Another difference is that in naive set theory one thinks globally: all sets are connected to each other in the sense that you can form unions, intersections and products of \emph{any} sets. In naive second order mathematics the basic thinking is local and structural: the structure of the natural numbers exists in countless isomorphic copies and none is more important than another. On any such structure we have sets, relations and functions. There is no obvious sense in which the subsets of one structure could be combined with the subsets of another. In set theory objects are conceived as elements of one domain, the domain of sets, and the domain itself is not a set. In second order mathematics there are numerous domains and we have second order mathematics on each domain separately. The domain is a set in the sense of second order mathematics, but there is no general theory of such domains, for none seems to be needed. We can start with as big domain as we wish.

The presence of many disconnected domains seems to contradict the idea of categoricity. How to relate a structure on one domain with structures on another domain? Categoricity can be made sense of by taking the union of the two domains as a new domain of individuals. Then the two structures would be parts of a unique structure which consists in effect of two structures. This is related to \emph{internal categoricity}, see section~\ref{internal}. Similarly, we can (in second order mathematics) define the natural numbers separately from the real numbers but we can also define them together so that the former are a substructure of the latter. Intuitively this has no effect on the natural numbers but in a finer analysis there is no guarantee of this. In the context of naive second order mathematics this potential non-uniqueness is not an issue, but it can arise in the context of  formalized second order logic. (See section~\ref{internal}.) In general, the locality of second order mathematics forces us to make more and more so-called \emph{large domain assumptions}, meaning that in order to investigate in a second order way a particular structure we have in mind we may have to conceive of it as part of a bigger structure, whether we want to or not.\footnote{Skolem  \cite{Skolem1923} (engl. trans. in \cite[p. 290]{MR209111}) criticizes ``domains".}
 In set theory this is not a problem because there is the global domain of all sets in the background all along. In second order logic there is no such global structure, except in second order set theory, which is a formalization of set theory using second order logic.

One way to compare the two naive foundational theories is to \emph{translate} them into each other. Indeed, naive second order mathematics can be routinely translated into set theory. The domain of individuals is chosen to form a set. Relations are interpreted as  set-theoretical relations on this set. Indeed,  this is quite generally considered the right way to think of second order mathematics to begin with. But of course, it is second order mathematics seen through the lens of naive set theory, and proponents of second order logic may see this as a betrayal. Likewise, there is no difficulty in building a set theory inside second order logic, as Zermelo did \cite{zbMATH02562682}.  
 Second order $\ZFC$, a kind of amalgam of second order logic and set theory, is an example of this. Proponents of (first order) set theory see here the problem, what is the meaning of the second order quantifiers?  One can also think of naive second order mathematics as a very weak naive set theory in which all sets are sets of individuals (urelements). By using relations to code families of relations and then further families of families of relations, on can build stronger and stronger naive set theory inside naive second order mathematics. This assumes that there are enough individuals to do the coding. For example, one $n+1$-ary relation $R$ and a set $A$ can code the family $\{S_a:a\in A\}$ of $n$-ary relations on a set $B$ via $$S_a(b_1,\ldots,b_n)\iff R(a,b_1,\ldots,b_n).$$ This presupposes a big enough set $A$. 

These two translations, second order mathematics into naive set theory on the one hand, and naive set theory into second order mathematics on the other hand, leave the basic question of comparing the two to each other on a neutral platform untouched. We either consider second order mathematics in the context of naive set theory or naive set theory in the context of second order mathematics. What is the more general context in which we can understand both equally and compare them? It would seem that ultimately set theory is the more general platform where both second order mathematics and naive set theory can be compared.
 
Admitting the difficulty of comparison, one still cannot help the feeling, while working with both,  that second order mathematics is just a fragment of naive set theory. This is because second order mathematics, whatever its formal language, is by its very nature a way to talk about mathematical objects and structures in their immediate neighborhood referring to subsets, relations and functions, while set theory is a way to talk about an entire universe which is rich enough to contain the mathematical objects and structures that mathematicians are interested in. On the other hand,  second order mathematics is usually considered a more economical theory of mathematics and, in comparison, naive set theory is sometimes considered wasteful, over-generating more objects than seem necessary. An example of this difference is that familiar mathematical structures such as the fields of real and complex numbers are defined in second order mathematics up to isomorphism only, while in naive set theory one defines these structures explicitly raising the question whether the explicit definition given is the right one. This difference is not a mathematical one but rather a question of taste, philosophy or aesthetics. Whatever is the right way to define familiar mathematical structures, second order mathematics or naive set theory, bigger differences emerge when we go deeper into mathematics, for example to the question whether all projective sets of reals are Lebesgue measurable. 

The formal language of second order logic (see Section~\ref{formalized}) is one which has variables for sets, relations and functions in addition to individuals. Otherwise the language of second order logic is like the formal language of first order logic. The formal language, on the other hand, of set theory has variables only for individuals, but these individuals are thought of as sets, sets of sets, sets of sets of sets, etc. The big difference between second order logic and set theory becomes apparent when we look at the axioms of both. The axioms of second order logic say that we can form new sets, relations and functions by means of definitions. The axioms of set theory, in contrast, stipulate that from any set we can form the power-set, union of its elements, definable subsets, and unions of definable families indexed by a set.  Roughly speaking, the axioms of set theory are more powerful than the axioms of second order logic, whence the latter gives the appearance of merely a weak fragment of the former.










\section{The story  of the  \emph{second order logic  --- set theory} relationship}\label{story}

Set theory was launched by Cantor \cite{zbMATH02716454} about 150 years ago.  A few years later  Frege 
\cite{zbMATH02709217} introduced a formal language which was essentially just second order logic. Thus second order logic and set theory emerged practically simultaneously.  Fifteen years after Cantor's paper \cite{zbMATH02716454}, Dedekind characterized the natural numbers with second order axioms \cite{zbMATH02693454}. A decade  later other mathematical structures, typially  the reals numbers, were characterized with second order axioms 
\cite{zbMATH02659709}, \cite{zbMATH02656983}, \cite{zbMATH02653775}. At about the same time, Zermelo proved the well-ordering principle with an argument that was given in essentially naive third order logic  \cite{zbMATH02652727}. Meanwhile paradoxes emerged in set theory and threatened its credibility. With the elimination of the paradoxes in mind,  Zermelo introduced his axiomatization which became the basis of our modern axiomatization of set theory \cite{zbMATH02639893}. Fifteen years later Fraenkel added the Replacement axiom schema and the axiom system $\ZFC$ was completed \cite{fraenkel}. It is an axiomatization of set theory in first order logic.

In the early foundational period\footnote{Represented by e.g. \cite{zbMATH02575336}.}, second order logic was considered a necessary and at the same time apparently innocuous extension of first order logic. However in 1930 G\"odel exposed the real nature of second order logic, when he pointed out that, unlike first order logic, second order logic is not (recursively) axiomatizable.\footnote{This was done verbally in the 1930 K\"onigsberg meeting and written up in {\em Erkenntnis 2}. See \cite{MR831941}, 1931a.} The reason for the non-axiomatizability is, as
  Gödel points out in his doctoral thesis, that ``For, if the unsolvability of some problem (in the domain
of real numbers, say) were proved, then, from the definition above, there would follow the existence of two non-isomorphic realizations of the axiom system for the real numbers, while on the other hand we can prove
the isomorphism of any two realizations" \cite[p. 62-62]{MR831941}. Gödel probably wants to emphasize here the difference between full second order logic and first order theories, without using this terminology. His completeness theorem applies only to the first order case. Curiously, this passage of the thesis was omitted from the published version \cite{MR1549799}. We may only guess what the reason was, but maybe he thought paying attention to the difference between first and second order logics would turn attention in the wrong direction, away from the main result of the thesis. On the other hand, this omission may have slowed down the recognition that second order logic is very different from first order logic.


Subsequent observations, involving concepts such as  L\"owenheim-Skolem Theorems, Hanf numbers and Compactness properties (\cite{MR71382}, \cite{MR204270}, \cite{MR295904}), led unequivocally to the conclusion that second order logic is a logic of a totally different nature than first order logic. It took some time before this difference in nature became duly appreciated.

 Kanamori argues: 
 \begin{quote}
``First and foremost, first-order logic is part and parcel of Gödel’s work
both in completeness and set theory. Zermelo proceeded in second-order terms
and ultimately did not take the linguistic turn, in that he did not develop
an uninterpreted formalism. Whereas the Skolem paradox much exercised
Zermelo, Gödel subsumed paradox into method by invoking Skolem’s analysis
to establish the continuum hypothesis in $L$. Gödel showed how first-order
definability can be formalized and used in a transfinite recursive construction
to establish striking new mathematical results. This significantly contributed
to a lasting ascendancy for first-order logic, which beyond its suﬃciency as
a logical framework for mathematics was newly seen to have considerable
operational eﬃcacy." \cite[p. 398]{MR2640544}
\end{quote}

Early axiomatizations of set theory, such as Bernays \cite{zbMATH03031148}, von Neumann \cite{zbMATH02575687}, Gödel \cite{zbMATH03099250},  were in fact second order, but their semantics was construed in terms of many-sorted logic, reminiscent of Henkin-models of SOL.
Zermelo's 1930 axiomatization $\ZFC$ of set theory in second order logic \cite{zbMATH02562682} was intended by its creator to be genuinely in (full) second order logic and this led to a famous  debate (see e.g. \cite{0998f894-5254-38b2-8237-b67692c3f8e7}, with Skolem who advocated a first order version, such as the present day $\ZFC$.

Fraenkel introduced what are now called permutation models in order to prove the unprovability of the Axiom of Choice from the other axioms of set theory, assuming we allow urelements \cite{zbMATH02600875}. A little later  Mostowski refined the method \cite{zbMATH03031149} and the resulting models became known as Fraenkel-Mostowski models. Mostowski worked in  second order logic. His 1938 doctoral thesis on the concept of finiteness was in  second order logic (or more exactly in type theory). Mostowski uses urelements, as is very natural in  second order logic. 

Axiomatic set theory arose from second and higher order logic, more generally from type theory when the typing was cast away.
Gödel writes in his 1933 unpublished lecture ``The present situation in the foundations
of mathematics":
\begin{quote}
It may seem as if another solution were afforded by the system of axioms
for the theory of aggregates, as presented by Zermelo, Frankel  and von Neumann; but it turns out that this system of axioms is nothing else but a
natural generalization of the theory of types, or rather, it is what becomes
of the theory of types if certain superfluous restrictions are removed. \cite[p. 34]{MR1332489}
\end{quote}

In the process of the transition from second order logic to set theory
 Quine \cite{bfc2ed3b-f285-3d09-b016-6dd131c63fed} showed that the simple theory of types in the sense of \cite{zbMATH02542491} can be embedded in set theory.
 Many who published in  second order logic moved to set theory after the Second World War, e.g. Mostowski \cite{zbMATH03031149}
  worked in second order logic but \cite{zbMATH03048504} and \cite{MR54547} and his subsequent papers are based on first order logic.
  
  On may speculate as to the reason why set theory started to dominate as the medium in which mathematical logic was presented. Was it because of philosophical considerations, perhaps a preference to present at least in the context of mathematical logic all of mathematics in one framework, even if practitioners of different parts of mathematics had no desire to represent their subject matter in set theoretic terms? Or was it because of the increasingly complicated notation inherent to second order logic and in its natural extension to type theory?

\section{Formalized second order mathematics}\label{formalized}

The formal language of second order mathematics, commonly known as second order logic, 
 was introduced by Frege and presented on a par with first order logic in the Hilbert-Ackermann classic \cite{zbMATH02575336}. In the much later but still classic monograph \cite{zbMATH03095463,zbMATH03122413} second order logic is presented properly as a highly non-trivial  extension of first order logic. 
 
Let us review for the sake of completeness what second order logic as a formal language actually is, as the details are not household material anymore today,  a century after the pioneering \cite{zbMATH02575336}. In short, starting with first order logic, we adopt new set-variables $X_n$, whose interpretation ranges over  subsets of the domain, new relation-variables $P^m_n$, whose interpretation ranges over  $n$-ary relations on the domain, and new function-variables $F^m_n$, whose interpretation ranges over total $n$-ary functions on the domain. New terms $F^m_n(t_1,\ldots, t_n)$ obtain inductively from the new function-variables, and new atomic formulas $X_n(t)$ as  well as  $P^m_n(t_1,\ldots,t_n)$ from the new set- and relation-variables, in combination with the new terms. We allow quantification over the new variables just as if they were individual variables of first order logic. 
A formula is said to be $\Sigma^1_0$, or equivalently $\Pi^1_0$, if it contains not second order quantifiers, $\Sigma^1_{n+1}$ if it starts with an existential second order quantifiers and the rest is $\Pi^1_n$, and 
$\Pi^1_{n+1}$ if it starts with an existential second order quantifiers and the rest is $\Sigma^1_n$.

As \emph{axioms} we take, following \cite{zbMATH02575336}, the usual propositional and quantifier axioms and the Comprehension Schemata, consisting of the formulas
\begin{equation}\label{comprehension}
\left\{
\begin{array}{l}
\exists X_n\forall x_1(X_n(x_1)\leftrightarrow\phi_1)\\
\exists P^m_n\forall x_1\ldots x_m(P^m_n(x_1,\ldots, x_m)\leftrightarrow\phi_2)\\
\exists F^m_n\forall x_1\ldots x_my(F^m_n(x_1,\ldots, x_m)=y\leftrightarrow\phi_3),
\end{array}\right.
\end{equation}
where $\phi_1$ and $\phi_2$ are  second order formulas with $x_1,\ldots,x_m$ but not $X_n$ (or $P^m_n$ in the case of $\phi_2$) free and $\phi_3$ is a second order formula with $x_1,\ldots,x_m,y$ but not $F^m_n$ free, and the universal closure of 
$$\forall x_1\ldots x_n\exists! y\phi$$ is made  an assumption\footnote{Hete ``$\exists !"$ is a shorthand for ``exists a unique".}. The Comprehension Schemata are a paradigm case of impredicativity in the sense that when e.g. a set (interpreting) $X_n$ is defined by (\ref{comprehension}),  the formula $\phi_1$ may have other bound set variables which intuitively range over all subsets of the domain, including the very set (interpreting) $X_n$ being defined. This has led to the introduction of limited versions of (\ref{comprehension}), most notably the $\Pi^1_1$-Comprehension Schema, where $\phi_1,\phi_2$ and $\phi_3$ are assumed to be $\Pi_1^1$-formulas. The role of the $\Pi^1_1$-Comprehension Schema is in second order logic comparable to the role of Kripke-Platek axioms KP (see e.g. \cite{MR424560}) in set theory.

Just as in first order logic, provability in second order logic can be equivalently formulated in terms of validity. The relevant concept of  model is the so-called Henkin model \cite{MR36188}, that is, a pair $(\mm,H_\mm)$, where $H_\mm$ is a set of subsets, $n$-ary relations and $n$-ary functions, $n<\omega$, on $M$ such that the schemata (\ref{comprehension}) hold. The truth of (\ref{comprehension}), or of any other second order sentence, in a Henkin model, is defined by letting the second order variables $X_n$, $P^m_n$, and  $F^m_n$ range over subsets, relations and functions, respectively,  that are elements of $H_\mm$. The corresponding concept in set theory is the concept of a transitive model of the $\ZFC$ axioms. Just as a transitive model probably misses some elements of the cumulative hierarchy up to its ordinal, a Henkin model $(\mm,H_\mm)$ probably misses some subsets, relations and functions on $M$. A second order analogue of the very particular, but at the same time very canonical, models $V_\alpha$ of $\ZFC$ are the \emph{full} Henkin models $(\mm,H_\mm)$, where $H_\mm$ is the set of \emph{all} subsets, $n$-ary relations and $n$-ary functions, $n<\omega$, on $M$. 

Second order logic with full models is called \emph{full second order logic}. As observed by Gödel (see above  Section~\ref{story}), there is no complete effective axiomatisation for full second order logic. This is a consequence of G\"odel's Incompleteness Theorem. The set theoretical analogue of full second order logic would be $\ZFC$ with semantics restricted to models of the form $V_\alpha$, $\alpha$ inaccessible, or equivalently, Zermelo's $\ZFC^2$ \cite{zbMATH02562682}, where Separation and Replacement Schemata are written in second order form and models are assumed to be full Henkin models. 

Since (\ref{comprehension}) is highly impredicative, it is not obvious how to construct Henkin models. However, we obtain Henkin models from models of set theory as follows: Suppose $K$ is a transitive model of $\ZFC$ and $\mm\in K$ is a structure. Let $H_\mm$ consist of all     subsets, $n$-ary relations and $n$-ary functions, $n<\omega$, of $M$ in $K$. Now $(\mm,H_\mm)$ is a Henkin model. The question of finding Henkin models can thus be reduced to the problem of finding models of set theory. For the latter problem there are the inner model and the forcing techniques. One may ask, should it not be easier to construct Henkin models than to construct models of set theory? It is easier in the sense that one can prove in set theory that Henkin models exist while one cannot prove that models of set theory exist. But fundamentally there is the same problem which arises from the inherent impredicativity of both second order logic and set theory. However, the Completeness Theorem of Henkin \cite{MR36188} provides an abundance of models, if we can assume consistency.
Thus there are consistent second order theories with a Henkin model but no full Henkin model, for example
\begin{equation}\label{incompact}
\{\phi_\oN(N,s,0)\}\cup\{N(c)\}\cup\{\neg c=0,\neg c=s(0),\neg c=s(s(0)),\ldots\},
\end{equation}
where $\phi_{\oN}(N,s,0)$ says that $s$ is a unary 1-1 function on $N$, $0\in N$, $0$ is not in the range of $s$, and every subset of $N$ which  contains $0$ and is closed under $s$, contains all of $N$. The non-full Henkin models, sometimes called \emph{non-standard} Henkin models of second order theories, are interesting in the same way as models, other than models of the form $V_\alpha$, of set theory.  It would be hard to study second order logic in full Henkin models only. In fact, the set of G\"odel numbers of second order sentences valid in all full models is $\Pi_2$-complete \cite{MR337481}, hence not $\Pi^m_n$ for any $m,n<\omega$. Respectively,  the set  of G\"odel numbers of first order sentences of set theory  valid in all $(V_\alpha,\in)$ is $\Pi_2$-complete\footnote{In contrast, the set  of G\"odel numbers of first order sentences of set theory  valid in all models of $\ZFC$ is $\Sigma^0_1$-complete.}.

A particularly important special case of formalized second order mathematics is second order number theory $Z_2$, an extension of Peano arithmetic to the real numbers, presented already in the monumental treatise  \cite{MR10509,zbMATH03334141} of Hilbert and Bernays. There is, up to isomorphism, only one full model of $Z_2$, so the interest is mainly in the Henkin models of $Z_2$. In this respect $Z_2$ resembles a (two-sorted) first order theory. A modern substantial account of $Z_2$ is in \cite{MR2517689}. For the complicated history of second order number theory, see \cite{MR3650982}. As a lot of theorems in ``ordinary" mathematics, such as calculus, can be represented and proved in second order number theory, we may consider $Z_2$ a partial foundation for mathematics, a second order  alternative to set theory. The consideration of subsystems of $Z_2$ can be used the delineate the strength of  such theorems. This has come to be called \emph{reverse mathematics}.
 Borel determinacy is an example, let alone stronger determinacy assumptions, such as projective determinacy $\PD$, of a statement expressible in the language of $Z_2$ and provable in $\ZFC$ but not provable in $Z_2$.

Another important special case of formalized second order mathematics is \emph{monadic} second order logic. Here ``monadic" means that the only second order variables that are allowed are the unary ones, and no function variables are allowed either. This kind of second order logic can be surprisingly weak on certain structures. Of course it is weak on structures with monadic predicates only. But it has been studied a lot on linear orders. This started with the result of B\"uchi to the effect that  monadic second order logic on $(\omega,<)$ is decidable  \cite{MR183636}. This was extended to decidability of two successor functions on the infinite binary tree \cite{MR246760}. On the other hand, ordinals bigger than $\omega$ have been considered in this context in \cite{Gurevich_Magidor_Shelah_1983} and \cite{doi:10.1142/9789812796554_0010}. The case of the decidability of monadic second order logic on $(\omega_2,<)$ has turned out to be particularly entangled with set theory.


Set theorists investigate arbitrary transitive models of set theory because there is a rich structure theory of such, involving inner models, ultrapowers and forcing extensions.  But what does such an investigation give us, considering that arbitrary transitive models are not the intended models $(V_\alpha,\in)$? The answer is  clear: investigating these arbitrary models tell us about the axioms, not about what is the case ``in reality", whatever that means. The models describe what the axioms really (can) tell us. They were instrumental in proving for example Borel Determinacy and they are instrumental in proving that Projective Determinacy follows from large cardinal assumptions. Thus, understanding arbitrary (``unintended") models gives us information about the intended universe $V$ of set theory, by helping us to see that some interesting things actually follow from the $\ZFC$-axioms {\em or their strengthening} by large cardinal assumptions. Since these things follow from  the axioms (plus large cardinals), we may conclude that they hold in $V$, subject to our faith in the large cardinals in question.

It is the same with second order logic. By working with Henkin models, we get an understanding of the axioms, such as the Comprehension Schema, and whatever new axioms we need, and this gives us important information about the 
full (i.e. intended) models. 

It is sometimes claimed that second order logic with Henkin models is just (many-sorted) first order logic. However, one should keep in mind that from a model theoretic point of view second order logic with Henkin models is much stronger than first order logic. Any second order theory with an infinite model is unstable.
This is because the Comprehension Schema allows one to build up complicated structures in the $H_\mm$ part of the model, even from the identity relation alone. So we should not confuse second order logic under Henkin semantics with first order logic. They are very different! If we drop  the Comprehension Schema assumption from Henkin models, the situation is different but then we have second order logic with no information about the second order variables and it is not surprising that the result is just first order logic.

One way to see that the Henkin semantics of second order logic is more reasonable than the full semantics, is the following: The concept of a Henkin model is absolute in set theory in this precise sense: If $(\mm,H_\mm)$ is a Henkin model of a second order sentence $\phi$ and $N$ is a transitive set such that $(\mm,H_\mm)\in N$ and $(N,\in)\models KP$, then $(N,\in)\models ``\mbox{$(\mm,H_\mm)$ is a Henkin model of $\phi$}"$. This is because the truth definition for Henkin models is essentially $\Sigma_0$: the first order variables range over $M$ and the second order variables range over $H_\mm$, both elements of $N$. This absoluteness means that second order logic with Henkin semantics is as independent of the ambient set theory as first order logic is. 

An alternative way to think of Henkin models is in terms of many-sorted first order logic. Instead of saying that $H_\mm$ is a set of subsets, $n$-ary relations and $n$-ary functions, $n<\omega$, on $M$, as in the concept of a Henkin model, one says that $H_\mm$ consists of 
\begin{enumerate}
\item A  set of abstract objects thought of as subsets and said to be of sort ``subset". 
\item A  set of abstract objects thought of as  $n$-ary relations and said to be of sort ``relation".
\item A  set of abstract objects thought of as  $n$-ary functions, for each $n<\omega$, on $M$ and said to be of sort ``function". 
\end{enumerate} Elements of $M$ are said to be of  sort ``individual". To connect the abstract objects, thought of as subsets, relations and functions, to real subsets, relations and functions, we have  binary ``epsilon"-relations between the elements of $M$ and the assumed abstract objects of the new sorts ``subset", ``relation" and ``function", for different arities. We denote these ``epsilon"-relations $E$, $E^r_n$ and $E^f_n$. The relations are supposed to satisfy the Axioms of Extensionality and the relevant Comprehension Schemata. For the ``epsilon"-relation relation $E$ the Axiom of Extensionality says
$$\forall x\forall y(\forall z (E(z,x)\leftrightarrow E(z,y))\to x=y),$$ where $x$ and $y$ are of the sort ``subset" and $z$ is of the sort ``individual", and similarly for the other new sorts. The Comprehension Schemata say
\begin{equation}\label{comprehension1}
\left\{
\begin{array}{l}
\exists x\forall z(E(z,x)\leftrightarrow\phi_1)\\
\exists x\forall z_1\ldots z_m(E^r_n(z_1,\ldots, z_m,x)\leftrightarrow\phi_2)\\
\exists x\forall z_1\ldots z_mu(E^f_n(z_1,\ldots, z_m,u,x)\leftrightarrow\phi_3),
\end{array}\right.
\end{equation}
where  $\phi_1$, $\phi_2$ and $\phi_3$ are formulas of the many-sorted (first order) language and $x$ is not free in them.

We can  ``Mostowski collapse" each abstract object $a$  to a real subset, relation or function $\pi(a)$ on $M$ as follows: We let $\pi(a)$ be
$$\begin{array}{ll}
\{z\in M: E(z,a)\}&\mbox{if $a$ is of sort ``subset"}\\
\{(z_1,\ldots,z_n)\in M^n: E^r_n(z_1,\ldots,z_n,a)\}&\mbox{if $a$ is of sort ``relation"}\\
\{(z_1,\ldots,z_n,z)\in M^n\times M: E^f_n(z_1,\ldots,z_n,z,a)\}&\mbox{if $a$ is of sort ``function"}.\\
\end{array}$$ Now $(M,\ran(\pi))$ is a Henkin model. It is easy to see that every Henkin model is of this form. This means that there is a perfect match between second order logic and many-sorted logic with the above sorts and the Axioms of Extensionality as well as the Comprehension Schemata  (\ref{comprehension1}).

The choice between the two---Henkin semantics and many-sorted logic---is a matter of taste. It is the question whether one talks about real subsets, relations and functions, or about things which are isomorphic to them, the Mostowski collapse providing the isomorphism. It is in the spirit of first order logic, and even more so in the spirit of second order logic, that we should not make a distinction between isomorphic objects. Logic should talk about structure and not about what things ``really" are. However, from the point of view of our basic framing---the comparison of second order logic and set theory---such a distinction seems most relevant. The fine line between Henkin semantics and many-sorted logic is a kind of dividing line, where second order ends and first order begins and this is certainly relevant from the point of view of comparing the two. Remarkably there is an isomorphism, namely the Mostowski collapse $\pi$, which builds a bridge between the two sides of the line. This is the root of the great difficulty in making strong statements about the difference between second order logic and set theory.   

One can go further and translate many-sorted logic into single-sorted logic by using unary predicates for the different sorts. This establishes second order logic in its many-sorted logic disguise as a particular first order theory, consisting of the translations of the Axioms of Extensionality as well as of the Comprehension Schemata. It should be noted that this does not mean that second order logic with Henkin semantics is the ``same" as first order logic. It only means that  second order logic with Henkin semantics is the ``same" as first order logic with certain set theory-like axioms.

Now that we have identified second order logic (with Henkin semantics) with a certain set theory-like first order theory, it becomes possible to compare the strengths of the former and $\ZFC$. Second order logic is like set theory with urelements, without the Axiom of Power-set and without the Axiom Schema of Replacement. It should be noted that second order logic is consistent with the axiom $\forall x\forall y(x=y)$. Thus it is much weaker than any reasonable set theory, unless we add  assumptions about the size of the universe. In \cite{MR2798269} such assumptions are called \emph{large domain assumptions}. Examples of such are the assumption of the infinity of the domain and the assumption of the domain having at least the size of the continuum. Without such assumptions we can prove the uniqueness of natural number and the real numbers, but we cannot prove their existence. When such large domain assumptions are translated into set theory, as above, they follow from the Axiom of Infinity and the Axiom of  Power-Set. Indeed, they can be seen as the prime reason for adopting such axioms among the axioms of set theory.

\section{Internal categoricity\label{internal}}

As we noted above, the categoricity of the concept of the natural numbers as well as of the ordered field of real numbers is a flagship feature of naive second order mathematics. This is how mathematics is built in the second order approach. At each step we stipulate what structure we are studying and since our stipulation is categorical, we ``know" what the structure is, or at least that is the second order mathematics' idea of knowing. 

When we formalize second order mathematics,  categoricity can be expressed as  a sentence which is provable from the axioms of second order logic. For natural numbers this means
the provability of \begin{equation}\label{cat-N}
(\phi_{\oN}(N,s,0)\wedge\phi_{\oN}(N',s',0'))\to\exists F\psi(F,N,s,0,N',s',0'),
\end{equation} where  $\phi_{\oN}(N,s,0)$ is as above in (\ref{incompact}) and $\psi(F,N,s,0,N',s',0')$ is the first order sentence which says that $F$ is an isomorphism between $(N,s,0)$ and $(N',s',0')$ \cite{MR3326591}. Note that proving (\ref{cat-N}) does not require defining the semantics of second order logic. 

The meaning of (\ref{cat-N}) is intuitively the same as our concept of the categoricity of the natural numbers. The provability of (\ref{cat-N}) from the axioms of second order logic  is a statement about finite sequences and is as unproblematic as any statement about the provability of a \emph{first order} sentence from \emph{first order} axioms. 


For the real numbers the categoricity means the provability of \begin{equation}\label{cat-R}
\begin{array}{l}
(\phi_{\oR}(R,+,\times,0,1,<)\wedge\phi_{\oR}(R',+',\times',0',1',<'))\to\\
\qquad\exists F\psi(F,R,+,\times,0,1,<,R',+',\times',0',1',<'),
\end{array}
\end{equation} where $\phi_{\oR}(R,+,\times,0,1,<)$ says that $+$, $\times$, $0$, $1$, and $<$ satisfy the axioms of completely ordered fields on $R$. Moreover, $\psi(F,R,+,\times,0,1,<,R',+',\times',0',1',<')$ is the first order sentence which says that $F$ is an isomorphism: $$F:(R,+,\times,0,1,<)\cong(R',+',\times',0',1',<').$$ This is the  textbook  definition, whenever one is given, as e.g. in \cite{MR1013117}, of the real numbers. Again, the intuitive meaning of (\ref{cat-R}) is  the same as our concept of the categoricity of the real numbers. The provability of (\ref{cat-R}) from the axioms of second order logic is, as in the case of (\ref{cat-N}), a statement about finite sequences and is as unproblematic as any statement about the provability of a first order sentence from first order axioms. 

The provability of (\ref{cat-N}) is called the \emph{internal categoricity} of  $\phi_{\oN}(N,s,0)$.
Respectively, the provability of (\ref{cat-R}) is called the \emph{internal categoricity} of the sentence $\phi_{\oR}(R,+,\times,0,1,<)$. For more on internal categoricity we refer to 
\cite{MR1867954}, \cite{MR2798269}, \cite{MR3326591},
\cite{MR3984974}, \cite{MR4329562}, and \cite{zbMATH07771690}. 

The word ``internal" stems from the following: The provability of a second order sentence, such as (\ref{cat-N}) or (\ref{cat-R}), means that the sentence is true in every Henkin model. Thus all models of (\ref{cat-N}) (respectively of (\ref{cat-R})) that are \emph{in} the same Henkin model, i.e. are \emph{internal} to the same Henkin model, are isomorphic. There are Henkin models with non-isomorphic natural numbers and real numbers when we cross from one Henkin model to another. More exactly, there are Henkin models  $(\mm,H_\mm)$ and $(\mm',H_{\mm'})$ of $\phi_\oN(N,s,0)$ such that 
$$(N,s,0)^{(\mm,H_\mm)}\ncong(N,s,0)^{(\mm',H_{\mm'})}.$$ This is a simple consequence of the Completeness Theorem proved by Henkin \cite{MR36188}. These two models $(\mm,H_\mm)$ and $(\mm',H_{\mm'})$ are not members of any single $H_\mn$. The Henkin model $(\mm,H_\mm)$ does not ``see" $(\mm',H_{\mm'})$ and vice versa.

Respectively, there are Henkin models  $(\mm,H_\mm)$ and $(\mm',H_{\mm'})$ of the sentence $\phi_\oR(R,+,\times,0,1,<)$ such that 
$$(R,+,\times,0,1,<)^{(\mm,H_\mm)}\ncong(R,+,\times,0,1,<)^{(\mm',H_{\mm'})}.$$ The so-called \emph{non-standard analysis\footnote{\cite{MR205854}}} is based on such a field $(R,+,\times,0,1,<)^{(\mm,H_{\mm})}$ which is non-isomorphic with the usual field of real numbers in the  sense that there are infinitesimals i.e. non-zero elements which are between $0$ and $1/n$ for each $n$.

Internal categoricity implies traditional categoricity because in full Henkin models every subset, relation and function is internal. Thus internal categoricity results imply the classical categoricity results. Even though the former are stronger than the latter, they are at the same time more convincing, because they do not depend on the definition of a semantics. They are finitist in the sense that they are results about finite sequences i.e. ultimately about natural numbers. 

There is another way to construe the classical categoricity results in a finitist way: They can be proved in $\ZFC$ and provability is a finitist property. However, internal categoricity results are results in the second order logic context. They do not depend on set theory. Provability in $\ZFC$ corresponds to internal categoricity when ``internal" is interpreted as internal to models of set theory rather than internal to Henkin models.

As Mostowski observed already in  \cite{MR66290}, the categoricity of second order definitions of familiar mathematical structures when interpreted in full Henkin models\footnote{Mostowski calls them `absolute' models.}, can be be recovered in arbitrary Henkin models if categoricity is redefined. Mostowski's redefinition agrees with  the received categoricity concept in full Henkin models. In non-full Henkin models Mostowski's categoricity means the ability to extend any isomorphism of the $H_\mm$-parts of the two models to an isomorphism of the first order parts.

\section{Model theory of fragments of second order logic}

Let us now consider full second order logic i.e. second order logic with its full Henkin models. This is an extension of first order logic which can be best understood and analyzed in set theory where fullness can be defined. If we did not use set theory to guarantee the fullness of the semantics, we could not prevent a plethora of non-standard (i.e. non-full) Henkin models creeping in as models of our second order sentences. In fact, first order logic, too,  can be best understood and analyzed in set theory as set theory provides a smooth theory of structures (models) for it. The appeal to set theory is much more relevant in the case of second order logic because the concept of a full Henkin model is non-absolute, being as it based on the power-set operation, unlike the concept of a  model for first order logic. Likewise, a carefully crafted second order sentence may have a model if and only if the Continuum Hypothesis is true. This would never happen for a first order sentence. A first order sentence  has a model if and only if there is no proof of contradiction from it. This property of the sentence is absolute in set theory. It is the same if we ask whether a second order sentence has a Henkin model. This is an absolute property of the sentence, thanks to Henkin's Completeness Theorem for second order logic with respect to Henkin models.

In the context of set theory, or with set theory as the metatheory, as it is also said, second order logic can be analyzed as an extension of first order logic without difficulty, even if its properties are very different. The goal of such an analysis is to extract as much information out of the fullness of the models as possible. As we have alluded to above, in the early days, i.e. before Gödel's results, the differences between first and second order logics were not recognized. Gradually they became generally known and recognized. Meanwhile a whole zoo of fragments of second order logic emerged and led to a general study of extensions of first order logic, such as infinitary languages and logics with generalized quantifiers, whether they are fragments of second order logic or not.

The first radical extension of first order logic, apart from second order logic, arose in Berkeley in the 1950s, where Tarski, Scott, followed by Henkin and Karp and motivated by problems in Boolean algebras, introduced  infinitary languages $L_{\kappa\lambda}$. The exact starting place was, according to Tarski, the logic seminar he had with Henkin in the fall of 1956. The results were published in \cite{MR99914} \cite{MR99915}, and later in \cite{MR166107}. Of course, infinitary languages are most certainly not literally fragments of second order logic, only in spirit, being as they are able to express non-first order properties such as finiteness, countability and well-foundedness. 
In an important further development in Berkeley,  Hanf announced results about compactness and Löwenheim-Skolem theorems for infinitary languages in \cite{Hanf1960}. He also introduced the concept of what came to be known as Hanf number of an infinitary logic. As Hanf  reports, the  currently well-known   more general definition of the Hanf number, which applies to any logic with a {\em set} of sentences, was suggested by Tarski. Although the questions about compactness and Hanf numbers of infinitary languages were in principle model theoretic questions, they have today  become household themes in set theory.

While infinitary languages were studied by Tarski's research group in Berkeley, in Poland Tarski's former student  Mostowski  introduced generalized quantifiers \cite{MR0089816}. Again, the goal was not to define fragments of second order logic but this was in effect what was going on. The most famous generalized quantifier introduced by Mostowski is undobtedly \begin{equation}\label{Q1}
\mm\models Q_1 x\phi(x)\iff |\{a\in M : \mm\models\phi(a)\}|\ge\aleph_1.
\end{equation} A complete axiomatization for the extension $L(Q_1)$ of first order logic by the obviously second order definable quantifier $Q_1$
was given by Keisler \cite{MR263616} after Vaught  \cite{MR166086} proved that such an axiomatization is in principle possible by showing that the set of Gödel numbers of valid sentences of $L(Q_1)$ is r.e.
The logic $L(Q_1)$ has also a Compactness Theorem in countable vocabularies as well as a Downward Löwenheim-Skolem Theorem down to $\aleph_1$, meaning that every model in a vocabulary of cardinality at most $\aleph_1$ has an $L(Q_1)$-elementary submodel of cardinality at most $\aleph_1$. 

If ``1" is replaced in (\ref{Q1}) by $``2"$, a different quantifier $Q_2$ obtains. The question of giving a complete axiomatization of $L(Q_2)$ (in $\ZFC$) is still open and properties of $L(Q_2)$ seem closely entangled with set theory. If cardinality is replaced in (\ref{Q1}) by cofinality, the following  ``cofinality quantifier", introduced by Shelah, obtains:
$$\mm\models Q^{cf}_\omega xy\phi(x,y,\vec{a})\iff$$
$$\{(c,d):\mm\models\phi(c,d,\vec{a})\}
\mbox{ is a linear order of cofinality $\omega$.}$$ This is a generalized quantifier in the more general sense introduced in \cite{MR244012}. Remarkably the fragment $L(Q^{cf}_\omega)$ of second order logic is axiomatizable and satisfies the Compactness Theorem in a vocabulary of any size  \cite{MR376334}. A different kind of modification of  (\ref{Q1}) is the so-called 
  Magidor-Malitz quantifier: 
$$\mm\models Q^{\mbox{\tiny MM}}_1 x_1,x_2\phi(x_1,x_2)\iff$$
$$\exists X\subseteq M(|X|\ge\aleph_1\wedge \forall a_1,a_2\in X:\mm\models\phi(a_1,a_2)).$$ The logic $L(Q^{\mbox{\tiny MM}}_1)$
is axiomatizable, assuming $\Diamond$ \cite{MR453484}, but consistently also non-axiomatizable \cite{MR1197203}. As with $L(Q_2)$, we see here an entanglement with set theory. Finally, so-called \emph{stationary logic} is receiving more and more attention. In this logic, denoted $L(aa)$, we have the new quantifier
$$\mm\models aa s \phi(s,\vec{t},\vec{a})\iff \{A\in \mathcal{P}_{\omega_1}(M)\ :\ (\mm,A)\models\phi(A,\vec{t},\vec{a})\}$$ $$\mbox{ contains a club of countable subsets of $M$.} $$ 
This sublogic of second order logic has a beautiful axiomatization in terms of the Fodor Lemma and it satisfies the Compactness Theorem in countable vocabularies \cite{MR376334},  \cite{MR486629}. The Downward Löwenheim-Skolem Theorem in the strong form (``Every model  in a vocabulary of cardinality at most $\aleph_1$ has an $L(aa)$-elementary submodel of cardinality at most $\aleph_1$.") is seriously entangled with stationary reflection in set theory, because it allows taking an arbitrary stationary subset $S$ of $\omega_2$ and finding an ordinal $\alpha<\omega_2$ of uncountable cofiality such that $S\cap\alpha$ is stationary. The claim is false if $V\!=\!L$ but consistently true, relative to the consistency of a supercompact cardinal (or in fact just $\mbox{PFA}^{++}$) \cite{MR506381}.

A general pattern in the above examples of extensions of first order logic is the following: we may obtain fragments of second order logic with model-theoretic properties resembling properties of first order logic but in almost all cases we face an entanglement with set theory in the sense that the model-theoretic property we are looking for holds in some models of set theory and fails in others. While the above extensions of first order logic were originally introduced in the context of model theory, they have gradually become part of set theory. The focus with these logics has moved from the desire to obtain an understanding of their model theory to the desire to understand set theory via them. These extended logics, sometimes also called \emph{strong logics}, play nowadays the role of suggesting new principles in set theory. 

We have discussed generalized quantifiers which have in one way or other similar properties as first order logic, even if this depends often on set theoretic principles. At the other end of the spectrum from strong logics are the ones that have similar properties as second order logic. An example is the game quantifier which is strong enough to say that a linear order is a well-order and on countable domains all $\Sigma^1_1$-properties. This, of course, kills axiomatizability and the Compactness Theorem but also pushes up the Hanf number into the range of large cardinals (if such exist) \cite{MR409188}. Another example is  the H\"artig quantifier  \cite{MR0189980} which is also strong enough to say, by means of extra predicates and extra sorts, that a linear order is a well-order  \cite{MR1136448}. If $V=L$, the H\"artig quantifier is as strong as full second order logic, in the sense that the two logics have the same $\Delta$-extension\footnote{See \cite{MR457146} for the definition and properties of the $\Delta$-extension.} but consistently the two are quite far from each other.

A good source for extensions of first order logic up to 1985 is  \cite{MR819531}. 
%
%
For the rich topic of set theoretic definability of logics we refer to \cite{MR819548}.

  



\section{Some uses of second order logic in set theory}

Second order logic and its fragments offer a criterion of definability in set theory. Here second order logic is taken with its full semantics. The most well-known example of this is descriptive set theory. The pointclasses ${\Sigma^1_n}$, ${\Pi^1_n}$ and ${\Delta^1_n}$, as well as their boldface versions $\boldsymbol{\Sigma}^1_n$, $\boldsymbol{\Pi}^1_n$ and $\boldsymbol{\Delta}^1_n$, are defined by reference to second order logic over the integers. They are nowadays an integral part of set theory. 

In the 1920s these pointclasses were already known but their investigation had hit a deadlock prompting the leading Russian researcher of such matters, Nikolai Luzin, to declare that a further progress is simply  impossible and will never happen. At this stage these pointclasses were investigated in topology, not set theory, which was in its infancy anyway. 

In the 1940s Kleene and Mostowski independently established the connection between hierarchies (e.g. Borel) studied in topology and quantifier prefix hierarchies arising in logic.  At about the same time, Gödel in his paper on the constructible hierarchy \cite{MR2514}, showed that the assumption $V\!=\!L$ leads to a breakthrough in the questions Luzin was so pessimistic about. For example, Gödel used $V\!=\!L$ to construct a 
$\boldsymbol{\Sigma}^1_2$ set of reals that is not Lebesgue measurable. This (or its negation) was exactly the kind of result Luzin was trying to prove. But Gödel used the axiom $V=L$, which does not immediately strike as a ``true" axiom, as it restricts sets to only those that can be constructed using Gödel's method. Why would all sets conform to such a restriction? In contrast, Solovay showed that consistently all $\boldsymbol{\Sigma}^1_1$ set of reals can be Lebesgue measurable \cite{MR265151}. Thus the question, which was originally a question about second order logic, had become one about the axioms of set theory. Eventually Martin and Steel \cite{MR959109}, in the culmination of a decades-long effort, proved that if the existence of infinitely many Woodin cardinals is assumed, all the above pointclasses have all the desired structural properties (Lebesgue-measurability, Baire property, perfect set property) investigated by Luzin and others in the 1920s. In particular all $\boldsymbol{\Sigma}^1_n$ sets of reals are Lebesgue measurable, whatever $n$ is. Thus a major result in set theory has culminated in structural properties of  sets of reals that are second order definable over the integers. The results were extended to $L(\oR)$ by Woodin \cite{MR959110} from a slightly stronger assumption. With $L(\oR)$ we move away from second order logic. But second order definability was never here a purpose in itself. It was just a stepping stone.
 
In a pivotal result about second order logic Magidor proved that it satisfies a strong form of the Downward Löwenheim-Skolem Theorem at a supercompact cardinal, and, moreover, a strong form of the Compactness Theorem at an extendible cardinal \cite{MR295904}. He proved: (1) If $\kappa$ is supercompact, $\phi$ is any second order sentence, and $M$ is any model, then there is a submodel $M'$ of $M$ of cardinality $<\kappa$ such that $\phi$ is true in $M'$, and (2) If $\kappa$ is extendible and $T$ is a second order theory in a vocabulary of any cardinality, and if every subtheory of $T$ of cardinality $<\kappa$ has a model, then $T$ has a model.   Magidor also showed that nothing less than a supercompact, respectively extendible, cardinal suffices for these results. This has been continued e.g. in \cite{MR2899689}. With this connection to supercompact and extendible cardinals second order logic had earned its keep in set theory. This also shows that if one is interested in second order logic and takes set theory seriously, one may obtain interesting results. It is not true that ``nothing" works in the model theory of second order logic, as it seemed up until the 1960s. 
 

 
 A recent development in set theory, pertaining to second order logic and its fragments, is the emergence of so-called \emph{extended constructibility}. Suppose $L^*$ is some fragment of second order logic. We use $C(L^*)$ to denote the modified Gödel's constructible hierarchy, obtained by replacing first order definability by definability in $L^*$. In an early result, Myhill and Scott proved that if $L^*$ is second order logic, then $C(L^*)$ is the model $\hod$ of hereditarily ordinal definable sets \cite{MR281603}. Thus, if $L^*$ is an arbitrary fragment of second order logic, then $$L\subseteq C(L^*)\subseteq \hod.$$  The inner model $C(L^*)$ has been particularly studied when $L^*$ is $L(Q^{cf}_\omega)$  \cite{MR4290501} and when $L^*$ is $L(aa)$ \cite{kmv2}. In both cases the theory of the inner model is absolute under forcing, assuming a proper class of Woodin cardinals. Thus it makes sense  to ask how do these inner models decide questions which are independent of $\ZFC$? These inner models are examples of how second order logic and its fragments may point to new directions of research in set theory. Second order logic and its fragments are weaker than set theory, to be sure, but in extended constructibility this weakness is turned into a strength. The weakness of these logics is taken advantage of in order to find order in the set theoretical universe.
 
 



\section{Summary}

In set theory, model theory, recursion theory, and philosophy of mathematics second order logic is usually, if not universally, understood to have the full Henkin semantics. This is in contrast with proof theory, where the semantics of second order logic, if any is provided,  is most naturally the completely general Henkin semantics. In set theory the full semantics of second order logic is assumed to be defined from the power-set operation of set theory. This renders second order logic the role of a restricted fragment of set theory essentially on the $\Sigma_2$-level of definability. This is a highly relevant level of definability where most if not all of ``ordinary" mathematics takes place. But the role of second order logic is here totally  auxiliary---important but auxiliary. It is auxiliary in the sense that it permits the formulation of results which are ultimately set-theoretic. If we look at second order logic from the proof theory perspective, it extends second order arithmetic from number theory to more general domains. Its proof-theoretic strength depends on the formulation of its basic concepts and on what assumptions of large domains are included, but it is  comparable in strength to Zermelo set theory. Whether we approach second order logic semantically or proof theoretically, it turns out to be a weaker cousin of set theory. If one likes the structuralist flavor of second order logic, then that is the way to go, and one can go a long way in analogy with set theory, but eventually one hits the notational complexities familiar from type theory.
If one accepts the universalist flavor of set theory, one can still use---and quite successfully---second order logic and its fragments, such as generalized quantifiers, as a restricted from of definability.
\bibliographystyle{plainurl}
\bibliography{SOLandST4Oxford}
\end{document}